\documentclass{article}
\usepackage{maple2e}
\usepackage{array}
\usepackage{amsmath}
\usepackage{graphicx}
\usepackage{url} 
\usepackage[top=.9in,bottom=.9in]{geometry}
 \numberwithin{equation}{section}
 
\DefineParaStyle{Maple Output}
\DefineCharStyle{2D Comment}
\DefineCharStyle{2D Math}
\DefineCharStyle{2D Output}
\begin{document}
\pagestyle{plain}
\begin{maplegroup}
\title{Notes on the Zeros of Riemann's Zeta Function}
\author{Michael S. Milgram \footnote{mike@geometrics-unlimited.com} \\
\it Consulting Physicist, Geometrics Unlimited, Ltd. \\
\it Box 1484, Deep River, Ont. Canada. K0J 1P0 }
\date{October 31, 2009}
\maketitle



\setcounter{page}{1}
\pagenumbering{arabic}


\begin{flushleft}

\textbf{Corrections added July 30, 2015.} 
\newline

1. The last few paragraphs of Section 2, beginning with the words  \textit{`There are two cases where (2.9) fails...'} cannot be valid.  This is demonstrated by a counterexample, and discussed more completely, in the appendix of a paper entitled \textit{Integral and Series Representations of Riemann's Zeta function, Dirichelet's Eta Function and a Medley of Related Results}, Journal of Mathematics, vol. 2013, Article ID 181724, 17 pages, (2013) {doi: 	10.1155/2013/181724}.  This shows that the original Abstract claim (below): `thereby establishing the truth of Riemann's hypothesis' is invalid. The remainder of this paper, particularly those sections dealing with the properties of $\zeta(1/2+i\rho)$ are valid and correct.
\newline

2. The right-hand side of equation (B.3) should be multipled by a factor -1/2.

\begin{abstract}
\vspace*{0.2cm}
The functional equation for Riemann's Zeta function is
studied, from which it is shown why all of the non-trivial, full-zeros of the Zeta function
$\zeta (s)$ will only occur on the critical line {$\sigma=1/2$} where {$s=\sigma+I\, \rho$}, 
thereby establishing the truth of Riemann's hypothesis.  Further, two
relatively simple transcendental equations are obtained; the numerical
solution of these equations locates all of the zeros of {$\zeta (s)$} on the critical line.

\end{abstract}


\section {Introduction}
The study of the non-trivial zeros of $\zeta(s)$ has been the
subject of myriad investigations over the years and is of ongoing interest in number theory. It has also recently received attention from the physics community \cite{physics}. 
Strangely, the results being presented here cannot be found in any of the summaries (e.g. \cite{Abram}, \cite{Edwards}, \cite{Bateman}, \cite{Ivic}, \cite{Patterson}, \cite{Sriv}, \cite{Titchmarsh})
or primary research articles\footnote{A short, only representative list!} (e.g. \cite{Kerimov}, \cite{Levinson}, \cite{Odlyzko}, \cite{Odlyzko2}) that I have consulted. Although it seems
inconceivable that they have escaped detection over the centuries, if
such is the case, a possible explanation is that the analysis involves
complicated manipulation of long expressions, a task best relegated to
computer algebra, and only in the last few years have
computer algebra codes reached a level of sophistication that allows such manipulation to
proceed. In any case, since these results (perhaps buried, more likely new) impart significant insight
into the nature and location of the zeros of $\zeta(s)$, I am taking the
opportunity to summarize here the results I have found.\newline 

On page 50 of Ivic's book (\cite{Ivic}),
it is written: "The functional equation for $\zeta(s)$ in a certain sense
characterizes it completely". Accepting the truth of that statement suggests that a
study of the functional equation should yield insight into the nature
of the zeros of {$\zeta(s)$}. That is the path taken here. 


\section {The functional equation inside the critical strip $0 \le \sigma \le 1 $}


The functional equation for {$\zeta(s)$} is well-known (e.g. \cite{Ivic}):

\begin{equation}
\zeta (1 - s)={\displaystyle \frac {2\,\Gamma (s)\,\mathrm{cos}(
{\displaystyle \frac {\pi \,s}{2}} )\,\zeta (s)}{(2\,\pi )^{s}}} 
\label{Zeta1}
\end{equation}

and the existence of the trivial zeros of $\zeta (-2n), n>1$ is immediately apparent due to the 
appearance of the cosine function on the right hand side.
With reference to Appendix A, where an index of notation will be
found, it is possible to break (1) into its real and imaginary parts,
giving the functional equation in an equivalent form:
\begin{equation}
{\widetilde \zeta _{I} (\sigma ,\rho ) }= - Q\,{\zeta _{I}(\sigma ,\rho ) } - P\,{\zeta _{R}(\sigma ,\rho ) }
\label{Zeta2a}
\end{equation}

\begin{equation}
{\widetilde \zeta _{R}(\sigma ,\rho ) }= - P\,{\zeta _{I}(\sigma ,\rho ) } + Q\,{\zeta _{R}(\sigma ,\rho ) }
\label{Zeta2b}
\end{equation}

where explicit expressions for the coefficient functions P and Q are
presented in Appendix B and reference to dependence on the independent (real) variables $\sigma$ and $\rho$, 
where {$s=\sigma+I\, \rho$} have been omitted. 

\vspace*{.2cm}
Instead of studying (\ref{Zeta1}), being a functional equation between complex variables 
and functions, consider the equivalent forms (\ref{Zeta2a}) and (\ref{Zeta2b}), which can be 
interpreted as the statement of a coupling that exists
among two independent functions ({$ \zeta _{R} (\sigma ,\rho ) $} and {$ \zeta _{I}(\sigma ,\rho )  $ }) and two dependent 
functions ({$ \widetilde \zeta _{R} (\sigma ,\rho ) $} and {$ \widetilde \zeta _{I}(\sigma ,\rho )  $ })
of two real variables $\sigma$ and $\rho$ . 
All quantities are real and this
is emphasized by writing {$ \zeta_R ( \sigma ,\rho ) $} to mean $ \Re ( \zeta ( \sigma + I \rho )\, ) $ 
and similarly for {$ \zeta_I ( \sigma ,\rho ) $} . The intent is to study (\ref{Zeta2a}) and (\ref{Zeta2b}) 
to determine if these two constraints can be used to specify a region(s) of the {$ (\sigma ,\rho ) $} plane 
(\,corresponding to the complex {$\it{s} $} plane\,) where full-zeros of {$ \zeta (s) $ } may possibly be found.   
"Half-zero" refers to points (or continuous regions) of the {$ (\sigma , \rho )$} plane where {$ \zeta _{R} (\sigma ,\rho )  =0 $} or {$ \zeta _{I}(\sigma ,\rho ) =0 $ } but not both; 
"full-zero"  refers to any of the set of points {$ (\sigma _{0}, \rho _{0})$} where {$ \zeta _{R} (\sigma _{0} ,\rho _{0}) =0 $} and {$ \zeta _I (\sigma _{0},\rho _{0} ) =0 $ } simultaneously.
Because {$ \zeta (s) $ } is known to be meromorphic (no branch cuts) \cite{Ivic}, the location of full-zeros of {$ \zeta (s) $ } must be isolated in the complex {\it{s}} plane, and 
this property will be reflected by a similar property of {$ \zeta _{R} (\sigma ,\rho ) $} and {$ \zeta _{I}(\sigma ,\rho )  $ } in the (real) $ ( \sigma , \rho )$ plane.

\vspace*{.2cm}
In the following, the intent is to search for zeros of {$ \zeta _R $ } and {$ \zeta _I $ } as a function of $\sigma $ 
with the variable $\rho $ being treated as a parameter ( {$ \rho = \rho _{p}$}). 
This corresponds to a search for full-zeros along horizontal lines of the $( \sigma , \rho )$ plane within the critical strip, 
graphically corresponding to that same strip in the complex plane $ s= \sigma + I \rho $.
\vspace*{.2cm}
Because P and Q have no singularities (poles)\footnote{except on the negative real axis which is outside the region of interest}, 
notice that if {$ \zeta _{R}(\sigma _{0} ,\rho _{0} )  =0 $} and {$ \zeta _{I} (\sigma _{0} ,\rho _{0} ) =0 $ }
at some point {$( \sigma_{0} , \rho_{0} )$} then (\ref{Zeta2a}) and (\ref{Zeta2b}) require that {$\widetilde \zeta _R (\sigma _{0} ,\rho _{0} ) = 0 $} 
and {$\widetilde \zeta _I (\sigma _{0} ,\rho _{0} ) =0$ }. So, any full-zero of  
{$ \zeta (s)$ } that lies in the range {$\sigma \leq 1/2 $} will be mirrored about the {$\sigma = 1/2 $} axis (the critical line) by a full-zero in 
the range {$\sigma \geq 1/2 $}, on the horizontal line $\rho = \rho _{p} $ . This property is well-known and does not necessarily hold true for half-zeros.

\vspace*{.2cm}
With this result in mind, a search constraint will be applied that imposes a necessary, but not sufficient condition for a zero of {$ \zeta (s)$ } to exist. 
That is, the functions {$ \zeta _R $ } and {$ \zeta _I $ } and their respective functions reflected about the critical line will be required to be equal (but not necessarily zero). 
A full-zero of {$ \zeta (s)$ } represents a special case of this more general condition. Specifically
\begin{equation}
{ \widetilde \zeta _{R}(\sigma ,\rho _{p})   =  {\zeta _{R}(\sigma ,\rho _{p} ) } }
\label{ZtReZR}
\end{equation}

  and 
\begin{equation}
{ \widetilde  \zeta _{I}(\sigma ,\rho _{p})  } = { \zeta _{I}(\sigma ,\rho _{p})  }   .
\label{ZtIeZI} 
 \end{equation}

Application of (\ref{ZtReZR}) and (\ref{ZtIeZI}) to (\ref{Zeta2a}) and (\ref{Zeta2b}) yields a set of
transcendental equations isolating correspondingly  special values of {$\sigma $} and {$\rho _{p} $} 
through the following constraints:

 \begin{equation}
{ \zeta _{R}(\sigma ,\rho _{p} )  = - \frac {(1+Q)}{P} \, \zeta _{I} (\sigma ,\rho _{p} )  }
\label{Zeta3a}
\end{equation}
 
 \begin{equation}
{ \zeta _{I} (\sigma ,\rho _{p})  = - \frac {(1-Q)}{P} \, \zeta _{R} (\sigma ,\rho _{p}) } ,
\label{Zeta3b}
\end{equation}
 giving a necessary condition on {$\sigma $} and {$\rho _{p} $} through the requirement  that 

\begin{equation}
{  P^{2} + Q^{2} = 1 }
\label{PQcondition}
\end{equation}
provided that 

\begin{equation}
{ \zeta _R (\sigma ,\rho _{p}) \neq 0  } \hspace*{.5cm} {\rm and/or} \hspace*{.5cm}  { \zeta _I(\sigma ,\rho _{p})  \neq 0  }.
\label{RIneq} 
\end{equation}

The cases corresponding to the failure of \ref{RIneq} will be discussed shortly.

\vspace*{.2cm}
For general values of {$\sigma $} and {$\rho $}, a lengthy calculation using (\ref{P}) and (\ref{Q}) shows that P and {\it Q} have the general property that

\begin{equation}
{  P^{2} + Q^{2} = (2 \pi )^{(1-2\, \sigma)}\,\frac {\rm{cosh}(\pi \rho )}{\pi } |\Gamma (\sigma + I \rho )|^{2}\:(1+\frac {\rm{cos}(\pi \sigma)}{\rm{cosh}(\pi \rho )}) }
\label{PQresult}
\end{equation}

from which (\ref{PQcondition}) imposes the following constraint on {$\sigma $} and {$\rho _{p} \rightarrow \rho $} after some rearrangement and the use of (\ref{gammaSq2}) :

\begin{equation}
{ (4\pi ^{2})^{{\sigma - \frac {1}{2}}} 
\frac{ \Gamma (\frac{1}{2} + I \rho) \Gamma (\frac{1}{2} - I \rho) }{ \Gamma (\sigma+ I \rho) \Gamma (\sigma - I \rho)  } - 1 
= \frac {\rm{cos}( \pi \sigma )}{\rm{cosh}( \pi \rho )} } 
\label{SRresult}
\end{equation}
for which the main solution is

\begin{equation}
\sigma = \frac {1}{2} ,\hspace*{.5cm} \rho \,\,\, {\rm arbitrary},
\label{sigeqhalf}   
\end{equation}
consistent with what Riemann famously hypothesized. See Appendix C where a second possibility is isolated and discarded.\newline

The converse is also true. That is, (\ref{sigeqhalf}) trivially implies the truth of (\ref{ZtReZR}) and (\ref{ZtIeZI}), but (\ref{PQcondition}) doesn't. 
But, with the exception of the case discussed in Appendix C, (\ref{sigeqhalf})  implies (\ref{PQcondition}) uniquely, so (\ref{sigeqhalf}) is a necessary and sufficient condition for all of 
(\ref{ZtReZR}) , (\ref{ZtIeZI}) and (\ref{PQcondition}), which themselves are prerequisites (necessary) for the presence of a zero of $\zeta (s)$. So, with the exception of the pathology discussed in Appendix C, (\ref{ZtReZR}) and (\ref{ZtIeZI}) can only occur, and hence a full-zero of {$ \zeta (s) $} can only be found,
when (\ref{sigeqhalf}) is satisfied, subject to (\ref{RIneq}), whose failure unfortunately corresponds to exactly those special values of {$\sigma $} and {$\rho $} of specific interest.

\vspace*{.2cm}
There are two cases where (\ref{RIneq}) fails - half-zeros and full-zeros. The case of half-zeros is easily dealt with, since it is clear that 
(\ref{Zeta2a}) and (\ref{Zeta2b}) are incompatible with (\ref{ZtReZR}) and (\ref{ZtIeZI}) at a half-zero unless {$P=0$} and {$Q=\pm 1$}, 
thereby satisfying (\ref{PQcondition}) spontaneously .
Thus there is no expectation that a half-zero will satisfy (\ref{Zeta3a}) and (\ref{Zeta3b}) in general, although the 
sieves (\ref{ZtReZR}) and/or (\ref{ZtIeZI}) may occasionally catch some half-zeros, so this case is a subset of the general result, and (\ref{sigeqhalf}) does not necessarily apply.

\vspace*{.2cm}
As noted before, all full-zeros of {$\zeta (s)$} are distinct, meaning that it is possible to expand {$\zeta (s)$} in a Taylor series in a neighbourhood of the full-zero.
Furthermore, the imaginary and real parts of a meromorphic function at a full-zero must be of the same degree, so for a full-zero of degree {\it m} in 
the neighbourhood of a solution to (\ref{Zeta3a}) and (\ref{Zeta3b}) where it happens that {$ \zeta _R (\sigma _{0} ,\rho _{p}) = 0 $}  and {$ \zeta _I(\sigma _{0},\rho _{p} )  = 0$}, one can write 

\begin{equation}
\zeta _{R} (\sigma ,\rho _{p}) = (\sigma - \sigma _{0})^{m} \zeta _{R}^{(m)}(\sigma _{0} ,\rho _{p}) /m! 
\label{tayR}
\end{equation}
and
\begin{equation}
\zeta _{I} (\sigma ,\rho _{p}) = (\sigma - \sigma _{0})^{m} \zeta _{I}^{(m)}(\sigma _{0} ,\rho _{p}) /m! 
\label{tayI}
\end{equation}

where 

\begin{equation}
\zeta _{R}^{(m)}(\sigma _{0} ,\rho _{p}) = {\frac {\partial ^{m} }{\partial \sigma ^{m}}} \zeta_{R} (\sigma ,\rho _{p}) |_{\sigma = \sigma _{0}}
\label{sigdefR}
\end{equation}
and
\begin{equation}
\zeta _{I}^{(m)}(\sigma _{0} ,\rho _{p}) = {\frac {\partial ^{m} }{\partial \sigma ^{m}}} \zeta_{I} (\sigma ,\rho _{p}) |_{\sigma = \sigma _{0}}
\label{sigdefI}
\end{equation}
It is emphasized that the partial derivative is taken with respect to {$\sigma$} because the search for a full-zero is
being conducted along a horizontal line in the ({$\sigma $},{$\rho $}) plane. Substitution of (\ref{tayR}) and (\ref{tayI}) into 
(\ref{Zeta3a}) and (\ref{Zeta3b}) yields (\ref{PQcondition}) and then (\ref{sigeqhalf}), the same result as before, except that the equivalent of (\ref{RIneq}) 
is always true, because {$\zeta _{R}^{(m)}(\sigma _{0} ,\rho _{p}) $} and {$\zeta _{I}^{(m)}(\sigma _{0} ,\rho _{p}) $} are non-zero by the definition of {\it{"a zero of degree m}"}.\newline

Thus, with the exception of the case discussed in Appendix C, (\ref{sigeqhalf}) is the only solution to a necessary condition for 
locating a full-zero of {$\zeta (s)$} in the finite ({$\sigma $},{$\rho $}) plane (and hence the finite complex {\it s} plane by extension), 
explaining why non-trivial, full-zeros of {$\zeta(s)$} have only ever been located
on the critical line (\ref{sigeqhalf}).

\section{On the critical line {$\sigma =1/2$}}

For the totality of this section and the next, the variable $\sigma = 1/2 $. With this understanding, the constraints (\ref{ZtReZR}) and (\ref{ZtIeZI}) 
reduce to an identity and (\ref{Zeta3a}) and (\ref{Zeta3b}) can conveniently be written in the  form

\begin{equation}
\zeta _{R} = \frac {N}{D_{R}} \zeta _{I}
\label{ZR}
\end{equation}
and
\begin{equation}
\zeta _{I} = \frac {N}{D_{I}} \zeta _{R}
\label{ZI}
\end{equation}

where expressions for N, $D_{R}$ and $D_{I}$ are given in Appendix B, yielding the further identities

\begin{equation}
N^{2} = D_{R}\,D_{I}
\label{nsq}
\end{equation}
and
\begin{equation}
D_{R} + D_{I} = 1\, ,
\label{DRpDI}
\end{equation}

 from which it is clear that 

\begin{equation}
0 \leq D_{R}\,{\rm{,}} \,D_{I} \leq 1.
\label{DRDIineq}
\end{equation}
since {$\it {D_{R}} $} and {$ \it {D_{I}} $} must have the same sign. (\ref{nsq}) demonstrates that
N shares the zeros of both {$\it {D_{R}} $} and {$\it {D_{I}} $}. Since both of the latter cannot
vanish simultaneously due to (\ref{DRpDI}), the zeros of {$\it {D_{R}} $} and {$\it {D_{I}} $} will locate
the half-zeros, but not the full-zeros, of {$\zeta (s) $} along the critical line, because if a zero of {$\it{N}$} 
carried one of the full-zeros of {$\zeta (s) $}, (\ref{ZR}) and (\ref{ZI}) show that the order of
the zeros of  {$\zeta _{R}$} and {$\zeta _{I}$} would be inconsistent. Specifically

\begin{eqnarray}
{D_{R} = 0 \Rightarrow D_{I}=1\, , \zeta _{I} = 0 \, ,  \zeta _{R} \neq 0 } \nonumber \\
{D_{I} = 0 \Rightarrow D_{R}=1\, , \zeta _{R} = 0 \, ,  \zeta _{I} \neq 0 } 
\label{DRDIZ}
\end{eqnarray}

The full-zeros of {$\zeta (s) $} for a zero of degree {$\it {m}$} 
are obtained by applying l'H\^{o}pital's rule of differentiation with respect to {$\rho $}. Any solution of

\begin{equation}
{\frac{\zeta _{R}}{\zeta _{I}} \Rightarrow \frac{\zeta _{R}^{(m)} } {\zeta _{I}^{(m)}}   } = \frac{N}{D_{R}}.
\label{zprime}  
\end{equation}

will thus isolate a potential full-zero of {$\zeta (s) $}, but as discussed previously, (\ref{zprime}) is only a necessary condition for achieving this task.
Thus a numerical solution does not guarantee that a full-zero has been found, although the set of all solutions will include all the full-zeros as a subset. Limited experimentation
(see Section 4) indicates that, at least for {$ m=1 $}, only the full-zeros of {$\zeta (s) $} are ever located by (\ref{zprime}); no solutions with {$ m=2 $} have been found. 

\section{Locating the Zeros}

The various functions introduced can be used to locate both the half- and full- zeros by numerically solving transcendental equations.
From (\ref{nsq}) and (\ref{DRDIZ}), all solutions of {$N^{2}=0$} will specify all the half-zeros of {$\zeta (s)$} on the critical line. 
In the notation of Appendix B,

\begin{equation}
 {C_{m}}\,\mathrm{cos}({\rho _{\pi }}) - {C_{p} 
}\,\mathrm{sin}({\rho _{\pi }})=0
\label{DZ}
\end{equation}
is a simple form of this constraint. Each successive solution with increasing values of {$\rho $} will locate successive 
half-zeros of {$\zeta _{R} $} and {$\zeta _{I} $} alternately, as illustrated in Figure (\ref{DZfig}). \newline

An interesting variant of (\ref{DZ}) arises by re-writing the terms explicitly, giving

\begin{equation}
\frac{\Gamma _{I}}{\Gamma _{R}} = \frac {\rm {tanh} (\pi \rho /2) + tan(\rho _{\pi })}   {1-\rm {tanh} (\pi \rho /2) \, tan(\rho _{\pi })},
\label{gratio}
\end{equation}  

and, to the extent that {$\rm {tanh}(\pi \rho /2) \approx 1 $}, (\ref{gratio}) can be inverted to read

\begin{equation}
\rm{tan}\it (\rho _{\pi } ) = \frac{ \Gamma _{I} / \Gamma _{R}-1 }  { \Gamma _{I}/ \Gamma _{R}+1 }.
\label{tanrp}
\end{equation}

If the first order Stirling's approximation (\cite{Abram}) for {$\rho \rightarrow \infty $} is applied to the ratio
 {$ {\Gamma _{I}/ \Gamma _{R}}$}, a simple form emerges:

\begin{equation}
\frac {\Gamma _{I}} {\Gamma _{R}} \approx -\rm{tan} (\rho - \rho _{L})
\label{approx1}
\end{equation}

which can replace the left-hand side of (\ref{gratio} ). Alternatively, (\ref{tanrp}) becomes

\begin{equation}
\mathrm{tan}({\rho _{\pi }}) \approx {\displaystyle \frac { - \mathrm{cos
}(2\,\rho ) + \mathrm{sin}(2\,{\rho _{L}})}{ - \mathrm{sin}(2\,
\rho ) + \mathrm{cos}(2\,{\rho _{L}})}} .
\label{taneq}
\end{equation}

These forms contain numerous poles and zeros and appear to have little numerical use, but may possibly be of use in deducing the spacing between zeros \cite{Odlyzko2},\cite{Saker}.\newline

The location of the full-zeros of {$\zeta (s) $} is specified indirectly in (\ref{zprime}). For simple zeros ({$\it{m}=1$}) the transcendental equation to be solved is

\begin{equation}
{\frac{\zeta _{I}^{\prime} } {\zeta _{R}^{\prime} }   } = -\frac{N}{D_{R}} ,
\label{z1p}  
\end{equation}

a more convenient form being
 
\begin{equation}
{D_{R}\,{\zeta _{I}^{\prime} } } + {N}\,{{\zeta _{R}^{\prime} }  =0 }.
\label{z2p}  
\end{equation}

from which the full-zeros\footnote{ and half-zeros belonging to $D _{R}=0$ )} can be found by standard numerical techniques (see figure (\ref{DZfig}). Although it may possibly be useful for numerical work, this
form is unsatisfying because it requires knowledge of the Zeta function derivatives, making it almost tautological. Unfortunately, a form for the full-zeros similar to
(\ref{DZ}), involving only the variable {$\rho $} and transcendental functions of that variable, eludes me.

\section{Summary}

The functional equation for {$\zeta (s)$} has been expressed in the form of a coupling between its real and imaginary components. It was shown that non-trivial, full zeros of {$\zeta (s)$}, if any exist, 
are only compatible with a solution to the functional coupling equations for special values of the underlying independent variable $"s"$\ . Two possible sets of values were located;
one of those regions has been explored by others and no zeros have ever been found. The remaining region consists of the critical line $s=1/2$. This establishes that Riemann's hypothesis is true. 
Additionally, two relatively simple transcendental equations were isolated, the zeros of which coincide with all the zeros of {$\zeta (s)$} on the critical line.

\section{Acknowledgements}

I am grateful to Vini Anghel, Dan Roubtsov and Bruce Winterbon for aid and discussion.


 
\appendix
\section{Appendix: Notation and identities}

The Riemann Zeta function {$\zeta (s)$} is written over the complex s plane:

\begin{equation}
s=\sigma + I \, \rho  \nonumber
\end{equation}
as
\begin{eqnarray}
\zeta (s) & = & \Re (\zeta (\sigma +I\,\rho)) + I\, \Im (\zeta (\sigma +I\,\rho)) \nonumber  \\
  & \equiv & \zeta _{R}(\sigma ,\rho ) + I\, \zeta _{I}(\sigma ,\rho ) \nonumber  \\
  & \rightarrow & \zeta _{R}+ I\, \zeta _{I} \:\: {\rm {when}} \:\: \sigma = \frac{1}{2},
\label{zetadef}
\end{eqnarray}
the latter for brevity. All variables in (\ref{zetadef}) are real. At reflected points, define
\begin{eqnarray}
\zeta (1-s) & = & \Re ( \zeta (1-\sigma -I\,\rho)) + I\, \Im ( \zeta (1-\sigma -I\,\rho)) \nonumber  \\
  & = & \zeta _{R}(1-\sigma ,\rho ) -I\, \zeta _{I}(1-\sigma ,\rho ) \nonumber  \\
 & \equiv & \widetilde \zeta _{R}(\sigma ,\rho ) -I\, \widetilde \zeta _{I}(\sigma ,\rho ) \nonumber  \\
  & \rightarrow  & \zeta _{R}- I\, \zeta _{I} \:\: {\rm {when}} \:\: \sigma = \frac{1}{2}.
\label{zetarefdef}
\end{eqnarray}

With reference to (\ref{zprime}), note that
\begin{eqnarray}
\zeta _{R}^{(m)} \equiv \frac {\partial ^{m} } {\partial \rho ^{m}} \zeta _{R} \\
\zeta _{I}^{(m)} \equiv \frac {\partial ^{m} } {\partial \rho ^{m}} \zeta _{I} \nonumber  
\label{partial}
\end{eqnarray} 
and, for $m=1$ 
\begin{eqnarray}
(\zeta _{R})^{\prime} \equiv \frac {\partial} {\partial \rho} (\zeta _{R}) = - \zeta ^{\prime} _{I} \\ 
(\zeta _{I})^{\prime} \equiv \frac {\partial} {\partial \rho} (\zeta _{I}) =  \zeta ^{\prime} _{R} \nonumber  
\label{zdiffs}
\end{eqnarray}

Similarly, the Gamma function is written
\begin{eqnarray}
\Gamma (s) & = & \Re (\Gamma (\sigma +I\,\rho)) + I\, \Im (\Gamma (\sigma +I\,\rho)) \nonumber  \\
  & \equiv & \Gamma _{R}(\sigma ,\rho ) + I\, \Gamma _{I}(\sigma ,\rho ) \nonumber  \\
  & \rightarrow & \Gamma _{R}+ I\, \Gamma _{I} \:\: {\rm {when}} \:\: \sigma = \frac{1}{2}.
\label{gammadef}
\end{eqnarray}

The following identities are noted \cite{Abram}

\begin{equation}
|\Gamma ( I \rho)|^{2} = \Gamma (I\rho )\,\Gamma (-I\rho) = \frac{\pi }{\rho \, \rm{sinh}(\pi \rho) }
\label{gammaSq0}
\end{equation}

\begin{equation}
|\Gamma (\frac{1}{2} +I\rho)|^{2} = (\Gamma _{R}+ I\, \Gamma _{I} )\,(\Gamma _{R}- I\, \Gamma _{I} )= \frac{\pi}{\rm{cosh}(\pi \rho) }
\label{gammaSq2}
\end{equation}

\begin{equation}
|\Gamma ( 1 +I \rho)|^{2} = \Gamma (1+I\rho )\,\Gamma (1-I\rho) = \frac{\pi \rho }{\rm{sinh}(\pi \rho) },
\label{gammaSq1}
\end{equation}

the symbols {$\rho _{\pi}$} and {$\rho _{L}$} are used to decrease the printed size of some formulae:

\begin{equation}
{\rho _{\pi} \equiv \rho \,\rm{log} (2 \pi ) }
\end{equation}

\begin{equation}
{\rho _{L} \equiv \frac {\rho }{2}\, \rm{log}(1/4+\rho ^{2})} ,
\end{equation}
and $m$ and $n$ are always positive integers.

\section{Appendix: Formulae}

In (\ref{Zeta2a}) and (\ref{Zeta2b}) the following functions are used

\begin{eqnarray}
\mapleinline{inert}{2d}{P3 :=
((-Gamma[R]*sin(rho[pi])+Gamma[I]*cos(rho[pi]))*cosh(1/2*Pi*rho)*cos(1
/2*Pi*sigma)+(-Gamma[R]*cos(rho[pi])-Gamma[I]*sin(rho[pi]))*sin(1/2*Pi
*sigma)*sinh(1/2*Pi*rho))*exp(-sigma*ln(2*Pi));}{%
\maplemultiline{
\mathit{P} = 2[( - {\Gamma _{R} (\sigma ,\rho ) }\,\mathrm{sin}({\rho _{\pi }})
 + {\Gamma _{I}(\sigma ,\rho ) }\,\mathrm{cos}({\rho _{\pi }}))\,\mathrm{cosh}(
{\displaystyle \frac {\pi \,\rho }{2}} )\,\mathrm{cos}(
{\displaystyle \frac {\pi \,\sigma }{2}} ) \\
\mbox{}\quad + ( - {\Gamma _{R}(\sigma ,\rho ) }\,\mathrm{cos}({\rho _{\pi }}) - {
\Gamma _{I}(\sigma ,\rho ) }\,\mathrm{sin}({\rho _{\pi }}))\,\mathrm{sin}(
{\displaystyle \frac {\pi \,\sigma }{2}} )\,\mathrm{sinh}(
{\displaystyle \frac {\pi \,\rho }{2}} )]\,e^{( - \sigma \,\mathrm{
ln}(2\,\pi ))} }
}
\label{P}
\end{eqnarray}

\begin{equation}
\mapleinline{inert}{2d}{Q3 :=
2*((Gamma[R]*cos(rho[pi])+Gamma[I]*sin(rho[pi]))*cosh(1/2*Pi*rho)*cos(
1/2*Pi*sigma)+(-Gamma[R]*sin(rho[pi])+Gamma[I]*cos(rho[pi]))*sinh(1/2*
Pi*rho)*sin(1/2*Pi*sigma))*exp(-sigma*ln(2*Pi));}{%
\maplemultiline{
\mathit{Q} = 2[({\Gamma _{R}(\sigma ,\rho ) }\,\mathrm{cos}({\rho _{\pi }}) + {
\Gamma _{I}(\sigma ,\rho ) }\,\mathrm{sin}({\rho _{\pi }}))\,\mathrm{cosh}(
{\displaystyle \frac {\pi \,\rho }{2}} )\,\mathrm{cos}(
{\displaystyle \frac {\pi \,\sigma }{2}} ) \\
\mbox{} \quad + ( - {\Gamma _{R}(\sigma ,\rho ) }\,\mathrm{sin}({\rho _{\pi }}) + {
\Gamma _{I}(\sigma ,\rho ) }\,\mathrm{cos}({\rho _{\pi }}))\,\mathrm{sinh}(
{\displaystyle \frac {\pi \,\rho }{2}} )\,\mathrm{sin}(
{\displaystyle \frac {\pi \,\sigma }{2}} )]\,e^{( - \sigma \,
\mathrm{ln}(2\,\pi ))} }
}
\label{Q}
\end{equation}

The following functions are introduced in (\ref{ZR}) and (\ref{ZI}):

\begin{equation}
\mathit{N} = {\displaystyle \frac {{C_{m}}\,\mathrm{cos}({\rho 
_{\pi }})}{\sqrt{\pi }}}  - {\displaystyle \frac {{C_{p}}\,
\mathrm{sin}({\rho _{\pi }})}{\sqrt{\pi }}}
\label{N} 
\end{equation}

\begin{equation}
\mathit{D_{R}} = {\displaystyle \frac {1}{2}}  - {\displaystyle 
\frac {1}{2}} \,{\displaystyle \frac {  {C_{p}}\,\mathrm{cos}({
\rho _{\pi }}) + {C_{m}}\,\mathrm{sin}({\rho _{\pi }})}{\sqrt{\pi
 }}} 
\label{DR}
\end{equation}

\begin{equation}
D_{I} = 1 - D_{R}    \nonumber
\end{equation}
where

\begin{eqnarray}
{C_{p}}& =& \mathrm{cosh}({\displaystyle \frac {\pi \,\rho }{2}} )\,{
\Gamma _{R}} + \mathrm{sinh}({\displaystyle \frac {\pi \,\rho }{2
}} )\,{\Gamma _{I}}\\
{C_{m}}& = & - \mathrm{sinh}({\displaystyle \frac {\pi \,\rho }{2}} )
\,{\Gamma _{R}} + \mathrm{cosh}({\displaystyle \frac {\pi \,\rho 
}{2}} )\,{\Gamma _{I}}
\label{CpCm}
\end{eqnarray}

\section{Appendix: Another solution?}

(\ref{sigeqhalf}) is the obvious solution to (\ref{SRresult}). Are there more? To answer this question
note that the magnitude of the right-hand side of (\ref{SRresult}) is strictly less than one, so any new solution with $\sigma \ne \frac{1}{2}$ 
must occur when the left-hand side is in that range. Consider {$L(\sigma ,\rho )$}, the left-hand side of 
(\ref{SRresult}) as a function of {$\rho $} at its endpoints {$\sigma =0 $} and {$\sigma =1$}. From (\ref{gammaSq0}) one gets

\begin{eqnarray}
L(0 ,\rho )& = & \frac {\rho }{2\pi } \rm{tanh}(\pi \rho) -1 \\
 & \rightarrow & \infty \quad \rm{as} \quad \rho \rightarrow \infty  \nonumber \\
 & \rightarrow & -1 \quad \rm{as} \quad \rho \rightarrow 0
\label{L0}
\end{eqnarray}

and from (\ref{gammaSq1}) one finds
\begin{eqnarray}
L(1 ,\rho )& = & \frac {2\pi }{\rho } \rm{tanh}(\pi \rho) -1 \\
 & \rightarrow & -1 \quad \rm{as} \quad \rho \rightarrow \infty  \nonumber \\
 & \rightarrow & 2\pi ^{2}-1 \quad \rm{as} \quad \rho \rightarrow 0
\label{L1}
\end{eqnarray}

Clearly {$L(\sigma ,\rho )$} changes sign for at least one value of {$\rho = \rho _{s} $} and $\sigma \ne\ \frac{1} {2}$, in the neighbourhood of which (\ref{SRresult}) could possibly be satisfied.
Numerically, {$ \rho _{s} = 6.283185307 $}; figure \ref{ldata} demonstrates that the slope of {$L(\sigma ,\rho )$} changes sign near {$\rho = \rho _{s} $} suggesting that
a numerical solution to (\ref{SRresult}) lies close by. To locate that neighbourhood precisely, consider

\begin{equation}
\frac{\partial}{\partial \sigma}L(\sigma , \rho ) = (L(\sigma ,\rho )+1) (4 \pi^{2}-2\Re (\psi (\sigma + I\rho)))\: .
\label{diffL}
\end{equation}

The sign of the left-hand side of (\ref{diffL}) will be determined by the factor

\begin{equation}
B(\sigma ,\rho ) = (4 \pi^{2}-2\Re (\psi (\sigma + I\rho))\: .
\label{bref}
\end{equation}

since the factor {$(L(\sigma ,\rho ) +1)$} is always positive. A change in the sign of {$B(\sigma ,\rho )$}
is consistent with the possibility of a numerical solution to (\ref{SRresult}). Figure \ref{bscan} shows that the sign of {$B(\sigma ,\rho )$} changes for
various values of {$\sigma$} and {$\rho$} near {$\rho _{s}$} with $0 \le \sigma \le 1,\, \sigma \ne \frac {1}{2}$, thereby isolating a second solution to (\ref{SRresult}), and a potential location to uncover 
a non-trivial, full-zero of {$\zeta (s)$} off the critical line. Others \cite{Odlyzko} have carefully searched this neighbourhood, and found no indication of such a zero. 
Since {$B(\sigma ,\rho )$} is monotonic with increasing(decreasing) values of {$\rho $}, there are no other possibilities. 
Thus (\ref{sigeqhalf}) defines the sole remaining range of possible solutions to (\ref{SRresult}).

\begin{figure}[p]
\begin{center}
\includegraphics[angle=0,scale=.8]{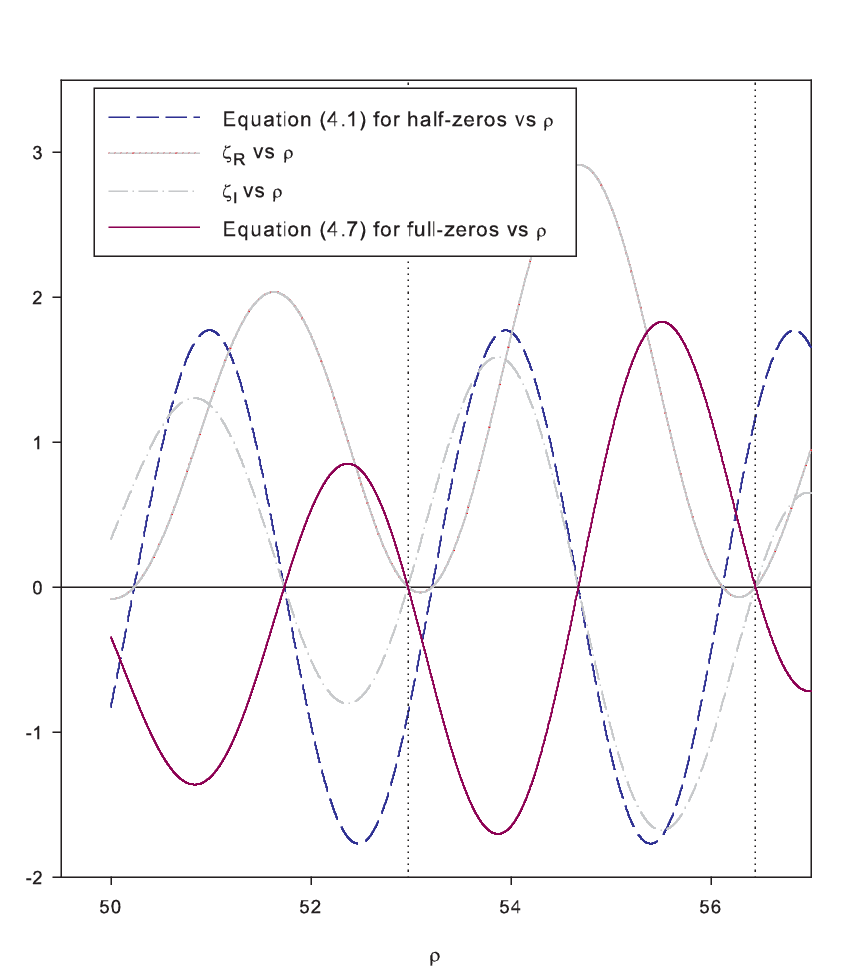} 
\end{center}
%
\caption {Numerical solution of (\ref{DZ}) and (\ref{z2p}) in the range {$50 \leq \rho \leq 57 $} showing conicidence with the half- and 
full- zeros of {$\zeta _{R} $} and {$\zeta _{I} $} respectively. The vertical dotted lines denote known full-zeros at {$\rho = 52.9703$}  
and {$\rho = 56.4462 $}.    }
\label{DZfig}
\end{figure}

\begin{figure}[p]
\begin{center}
\includegraphics[angle=0,scale=0.9]{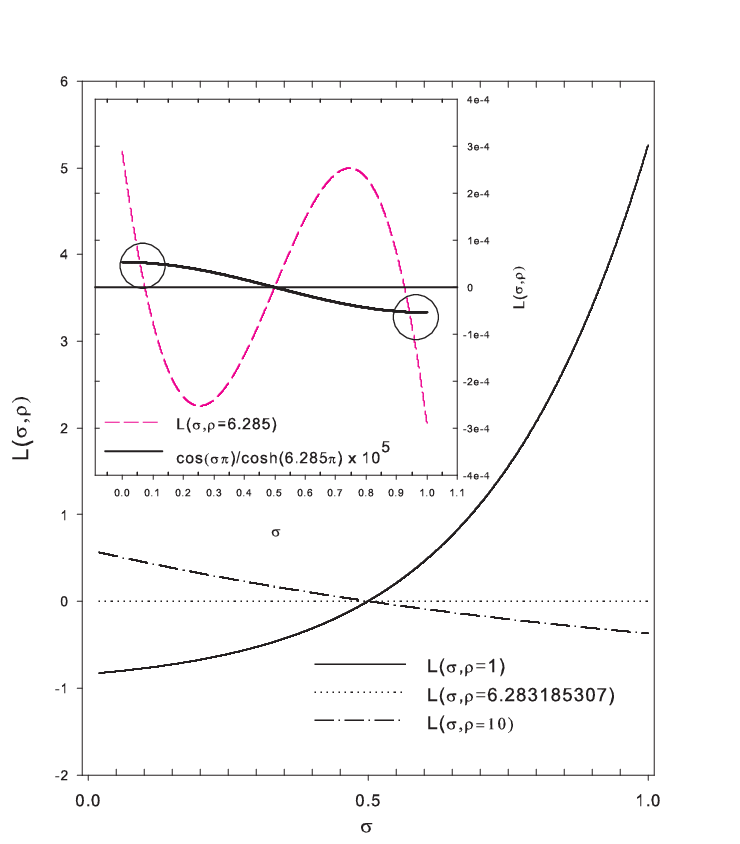} 
\end{center}
\caption{A plot of {$L(\sigma ,\rho )$} at three different values of {$\rho $}, bracketing {$\rho _{s} $}.
The inset contains a $10^{5}$ magnification of the right-hand side of (\ref{SRresult}), and the circles 
indicate intersection points of the two curves with {$\sigma \neq 1/2 $}, 
yielding a numerical solution to (\ref{SRresult}) and the location of a potential zero of {$\zeta (s)$} off the critical line. 
}
\label{ldata}
\end{figure}

 
\begin{figure}[p]
\begin{center}
\includegraphics[angle=0,scale=.6]{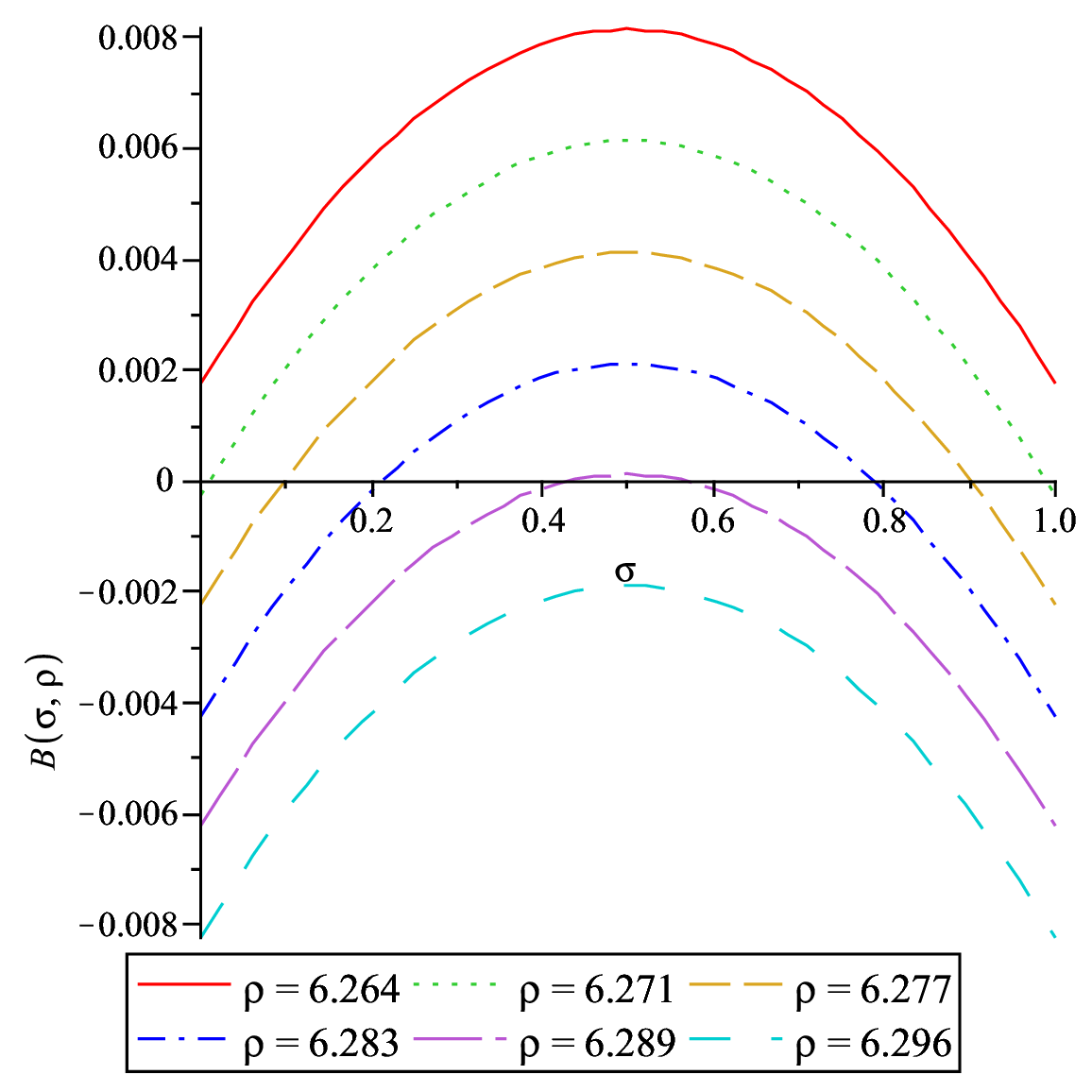} 
\end{center}
\caption{A parametric scan of the function B({$\sigma ,\rho $}) near {$\rho _{s} $} as a function of {$0 \leq \sigma \leq 1$}
}
\label{bscan}
\end{figure}

\end{flushleft}

\end{maplegroup}
\end{document}